\documentclass[runningheads]{llncs}
\usepackage[table,xcdraw]{xcolor}
\usepackage{graphicx}
\usepackage{amsmath}
\usepackage{amssymb}
\usepackage{subcaption}
\begin{document}
\title{Partial differential equation solver based on optimization methods}
%
%
\author{Alexander Hvatov}
\authorrunning{A. Hvatov}
%
\institute{ITMO University, Kronverksky pr. 49, Saint-Petersburg, Russia\\
\email{alex\_hvatov@itmo.ru}\\
}
\maketitle              
\begin{abstract}
The numerical solution methods for partial differential equation (PDE) solution allow obtaining a discrete field that converges towards the solution if the method is applied to the correct problem. Nevertheless, the numerical methods usually have the restricted class of the equations, on which the convergence is proved. Only a small amount of "cheap and dirty" numerical methods converge on a wide class of equations with the lower approximation order price. In the article, we present a method that uses an optimization algorithm to obtain a solution that could be used as the initial guess for the wide class of equations.

\end{abstract}
\section{Introduction}
\label{sec:intro}

Partial differential equations (PDE) solution is a traditional topic in the mathematical physics and applications \cite{morton2005numerical}. A wide variety of methods starting from finite-difference schemes through finite element method to modern spectral-like analytical methods, are established. However, most of the methods are aimed to solve the given equation. That means that the operator properties and possible boundary condition types are a priori known.

Modern data-driven methods \cite{maslyaev2020data} provide data-driven discovery of the partial differential equations to obtain possibly interpretable models. Equation discovery may result in an equation that, in the general case, does not have the known properties. An expert may solve the obtained equation. However, it is required to spend the expert's time every time the equation is discovered. Additionally, for the given data, we may use some shortcuts and obtain the solution automatically. However, the obtained solution is particular for the train observational field is of less interest.

 The discovered equation may be used further, for example, as part of the ensemble. Thus, it must have the solution in general form with other sets of boundary and initial conditions. Since the resulting equation properties and boundary conditions types are unknown beforehand, the specific PDE solution methods are not viable. 

Existing methods that may solve arbitrary equations either significantly reduce the class of the equations, as spectral methods, or require retraining of the neural network, which is, above all, time-consuming. In the article, we aim to obtain a "cheap and dirty" method of PDE solution, which possibly combines the broad operator class and high level of automatization of the process. It means that we, without the help of the expert, try to obtain the field that approximates the equation's solution. 

The paper is organized as follows: Sec.~\ref{sec:literature_rewiev} contains a brief review of the existing PDE solution methods, Sec.~\ref{sec:problem_statement} contains the definitions and algorithm description used in the article, Sec.~\ref{sec:numerical_exps} contains the application of the given algorithm to particular PDEs, Sec.~\ref{sec:conclusions} outlines the paper and proposes the directions for the future work.

\section{Related work}
\label{sec:literature_rewiev}

The classical finite-difference \cite{thomas2013numerical} and finite-element method \cite{solin2005partial} (FEM) have established area of applicability. For example, FEM is widely used to solve elliptic equations occurring in different areas, for example, mechanics.

Without a doubt, decades of development made FEM the fast method to solve known physics, and mechanics-related problems \cite{pavlovic2020geometry}. However, there is no possibility to apply finite-difference and finite element methods to arbitrary equations. Finite-difference methods could be applied to the linear equations. However, it is required to derive schemes to every order of the derivative that occurs.

Spectral methods for PDE solution are the most modern analytical and numerical methods \cite{burns2020dedalus}. However, their application to an arbitrary problem is restricted by the linear differential operators and boundary conditions. Moreover, it uses automatic differentiation on the polynomial decomposition series that also restricts the solutions' class.

Arising neural differential operators methods are slightly dependent on a training dataset and require to learn neural network every time a new problem arises \cite{li2020neural}. However, the recent research shows that combined with the transition to the spectral domain may be promising \cite{li2020fourier} even though it also restricts the applicability to the linear methods if applied to the Fourier specter directly.

\section{Problem statement}
\label{sec:problem_statement}

We solve the boundary PDE problem defined on a subdomain $(x,t) \in \Omega \subset 	\mathbb{R} ^2$ with a boundary $\partial \Omega$ in form Eq.~\ref{eq:PDE_problem}. We emphasize that the approach will work for higher dimensions. However, in illustrative matters, we show two dimensional (in examples below, we use a single space and single time dimension) equation.

\begin{equation}
\begin{array}{cc}
Lu=f\\
Bu=g 
\end{array}
\label{eq:PDE_problem}
\end{equation}

In Eq.~\ref{eq:PDE_problem} we assume that the differential operator $L$ and the boundary operator $B$ and the arbitrary functions $f,g$ are defined such that the boundary problem is correct. 

Most of the numerical methods assume that the solution field is found in a discrete subspace in form of the mesh function, it means that:

\begin{equation}
\begin{array}{cc}
\bar{u}=\{u(x^{(i)},t^{(i)}), i=1,2...,n\}\\
\forall i \, (x^{(i)},t^{(i)}) \in \Omega
\end{array}
\label{eq:field_discretization}
\end{equation}

Without loss of generality, we assume that the field discretization $X\allowbreak=\allowbreak\{x^{(i)}\allowbreak,\allowbreak t^{(i)}\}\subset\Omega$ is fixed during the process of PDE solution. For the experiments, we use a uniform mesh. However, the discretization for the method described below could be chosen arbitrarily.

Classical numerical methods assume that the values are connected using either finite difference schemes or variational principle as in the finite element method.

We formulate a minimization problem to find the solution field as Eq.~\ref{eq:prec_algorithm_formulation}.

\begin{equation}
    \min \limits_{\bar{u}} || L\bar{u}-f||+\lambda ||B \bar{u}-g||
    \label{eq:prec_algorithm_formulation}
\end{equation}

In Eq.~\ref{eq:prec_algorithm_formulation} $L$ is assumed to be the ``precise'' operator that gives the exact value of the derivative at the mesh points. We note that Eq.~\ref{eq:prec_algorithm_formulation}, $\lambda$ is an arbitrary chosen constant, which, if the boundary conditions are correctly defined, does not influence a result. In this case, there is no doubt that the optimum will be the solution to the differential operator. 

In practice, differential and boundary operators are also the approximation of the derivative that has an error, and the minimization algorithm is the numerical algorithm that has its error for different optimization problems. Therefore, the final problem that is solved in the article is formulated as Eq.~\ref{eq:approx_algorithm_formulation}.

\begin{equation}
    \min \limits_{\bar{u}} \sum  \limits_{\bar{x} \in X} ( \bar{L}\bar{u}-f(\bar{x}))^2+\lambda \sum  \limits_{\bar{x} \in X} ( \bar{B}\bar{u}-g(\bar{x}))^2
    \label{eq:approx_algorithm_formulation}
\end{equation}

In Eq.~\ref{eq:approx_algorithm_formulation} $\bar{L}$ and $\bar{B}$ are the approximate differential and boundary operators (meaning that the derivatives are replaced with the approximations), $\bar{x}$ are grid points and $\bar{u}$ are taken accordingly the given grid point.

\section{Numerical experiments}
\label{sec:numerical_exps}

It is necessary to prove two of the three following properties: convergence, stability, and approximation to prove the correctness of every numerical algorithm. In the article, we do not pursue the goal to prove that for an arbitrary operator $L$. This section shows several numerical experiments that can be used as the starting point and proof-of-concept.

\subsection{Practial realization}

As the derivative approximation, we use the finite-difference scheme of the second-order (approximation order $O(h^2)$, where $h$ is the uniform grid step in the discretization of the given dimension). We use both forward and backward for boundaries in the form Eq.~\ref{eq:forward_scheme_first}.

\begin{equation}
   \begin{array}{cc}
    u'_f(x) = \frac{u(x+h)-u(x)}{h} \\
    u'_b(x) = \frac{u(x)-u(x-h)}{h}
   \end{array}
   \label{eq:forward_scheme_first}
\end{equation}

For the interior points we use scheme Eq.~\ref{eq:center_scheme_first} as more stable.

\begin{equation}
u'_c(x) =\frac{1}{2}(u'_f(x) + u'_b(x)) =\frac{u(x+h)-u(x-h)}{2h}
   \label{eq:center_scheme_first}
\end{equation}

For the higher-dimensional derivatives, we apply the same scheme several times. Even though such an approach leads to a high approximation error, we intend to use proof-of-concept, which works even with this setup. Without a doubt, direct usage of the higher derivative order schemes will lead to a better result.

After the approximate differential operator is defined, we make the procedure that allows us to apply the arbitrary differential operator and boundary operator to the arbitrary field. Basically, we encode every operator with the axis and number of Eq.~\ref{eq:forward_scheme_first}-Eq.~\ref{eq:center_scheme_first} application. Also, we must encode the coefficient before every differential term and the power of the term.

Starting the arbitrary field, we use the optimization algorithm to minimize the difference between the applied operator to the field and function $f$ over all discretization points. Additionally, we introduce the difference between the applied boundary operator and function $g$. The last complimentary allows us to solve homogeneous equations non-trivially.

\begin{figure}[h!]
 \centering
 \includegraphics[width=0.95\linewidth]{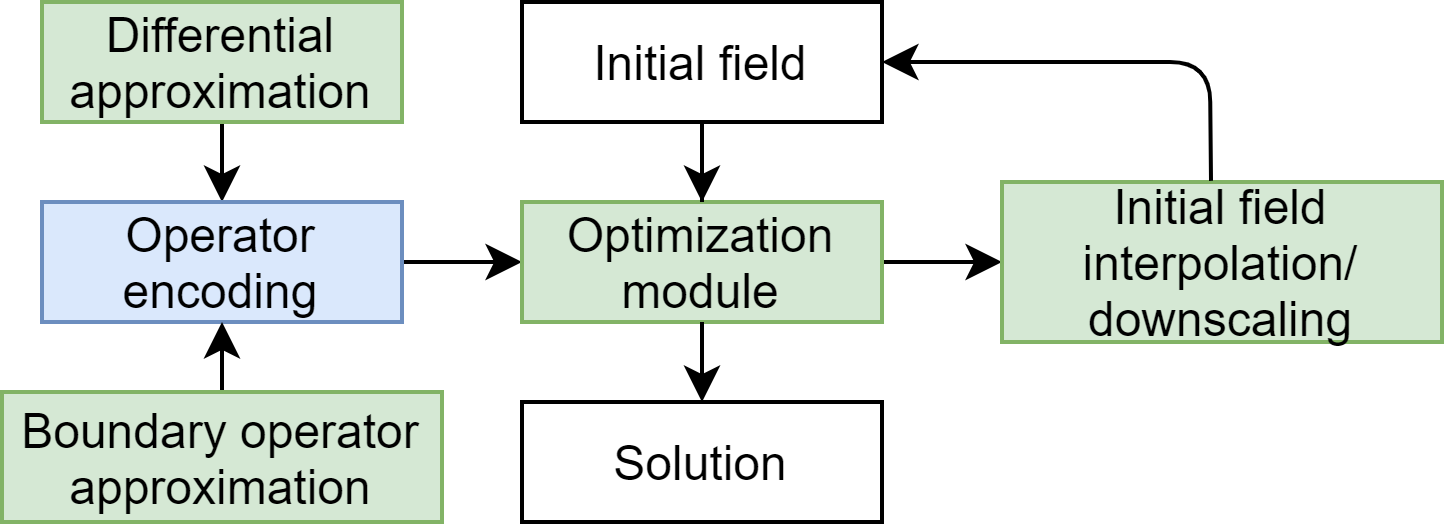}
\caption{Module structure of the solver algorithm green parts are replaceable.}
\label{fig:modules}
\end{figure}

We propose the module structure shown in Fig.~\ref{fig:modules} of the resulting solver.

Below we discuss how the module's replacement affects the resulting solutions and the time of the optimization.

\subsection{Convergence and stability}

This subsection shows a convergence of the algorithm for two examples: wave equation of the hyperbolic type and heat transfer equation of the parabolic type. 

\subsubsection*{Wave equation}

We try to assess the convergence of the algorithm the solution of the wave equation with boundary conditions in form Eq.~\ref{eq:WE_formulation}

\begin{equation}
    \begin{array}{cc}
         \frac{\partial^2 u(x,t)}{\partial t^2}-\frac{1}{4} \frac{\partial^2 u(x,t)}{\partial x^2}=0 \\
         u(0,t)=u(1,t)=0 \\
         u(x,0)=u(x,1)=\sin(\pi x) \\
         (x,t) \in [0,1] \times [0,1]=\Omega
    \end{array}
    \label{eq:WE_formulation}
\end{equation}

We use the formulation Eq.\ref{eq:prec_algorithm_formulation} to obtain the solution of the equation for 30 runs for consequently increasing the number of points in discretization from $10 \times 10$ points in $\Omega$ (since the mesh is assumed uniform, it is equal to $h=1/(N-1)=1/9$, where $N$ is the number of points for time and space dimensions, i.e., we take 10 points in the range $[0,1]$ including boundaries) to $50 \times 50$ points with the step of 5 points. We initialize optimization with the random field from the uniform distribution on a field value range $[0,1]$ for each run. 

\begin{figure}[h!]
\centering
\begin{subfigure}{.5\textwidth}
  \centering
  \includegraphics[width=0.95\linewidth]{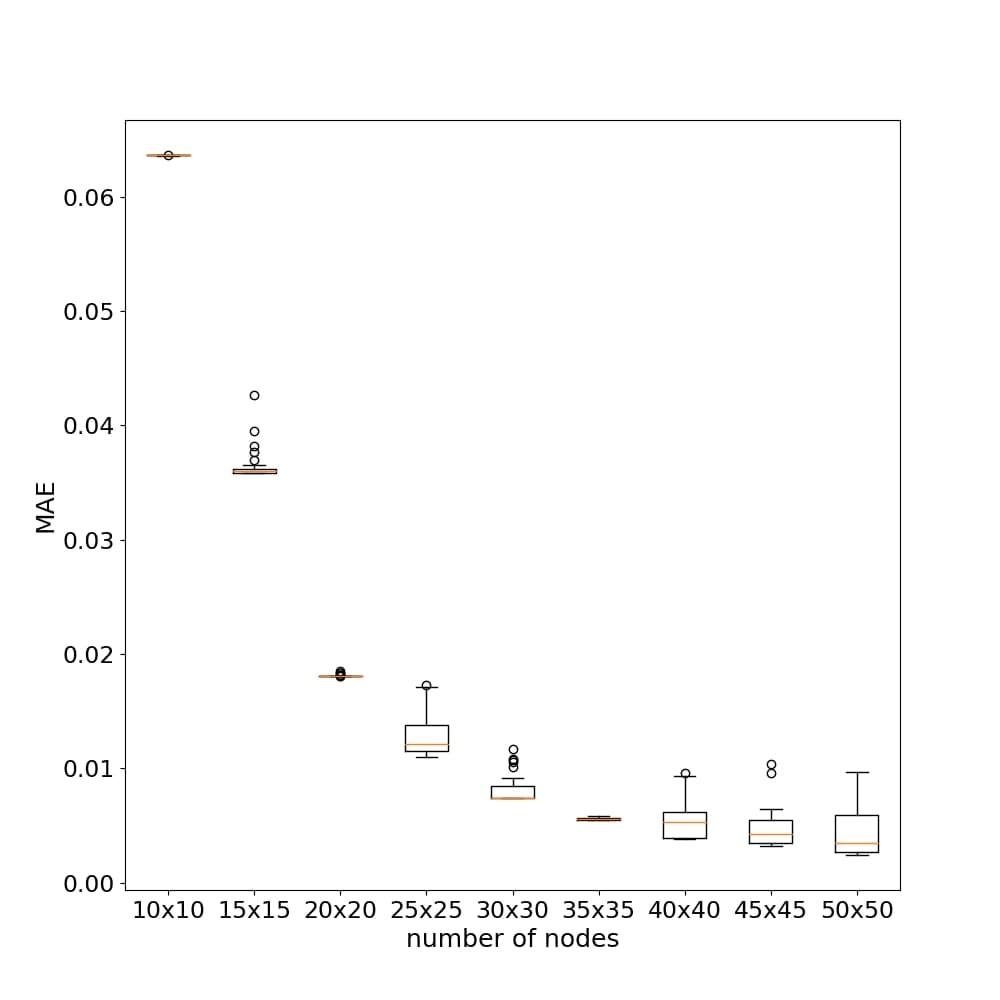}
  \caption{MAE}
\end{subfigure}%
\begin{subfigure}{.5\textwidth}
  \centering
  \includegraphics[width=0.95\linewidth]{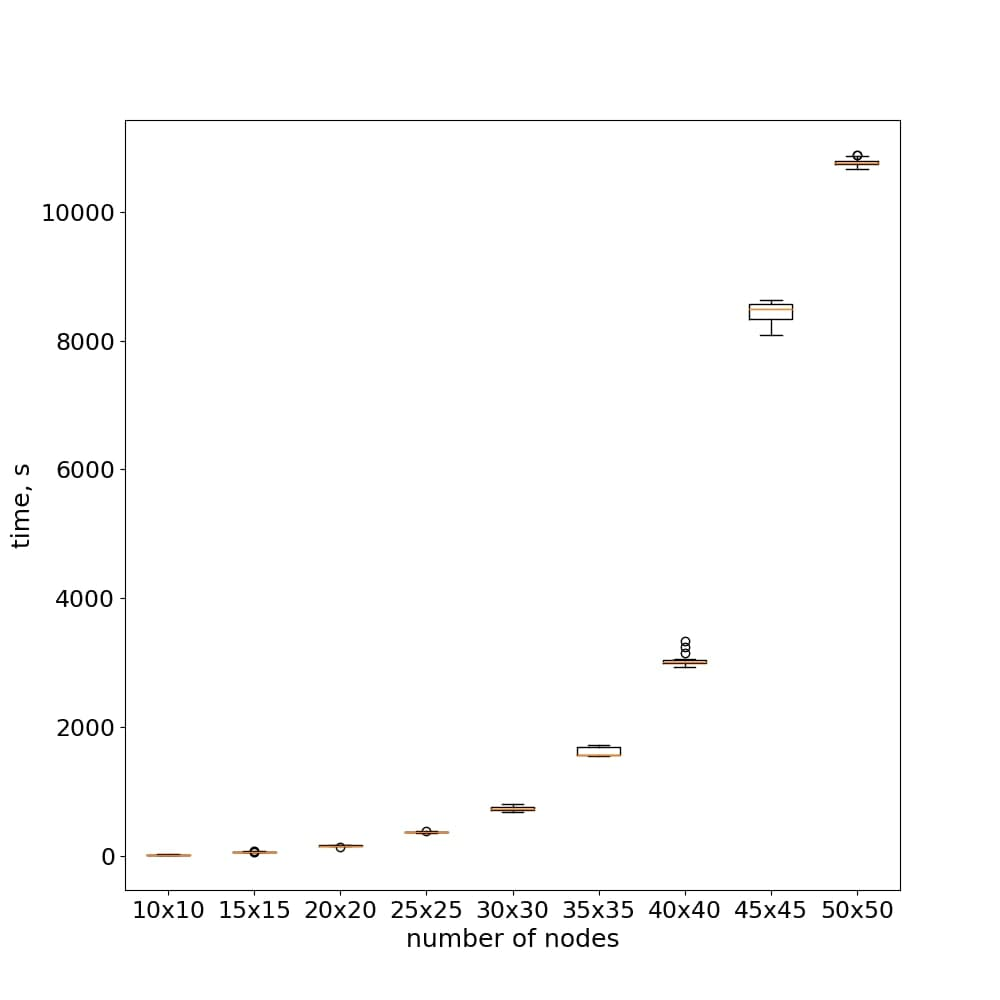}
  \caption{Optimization time}
\end{subfigure}
\caption{Convergency and time complexity for the wave equation (Eq.~\ref{eq:WE_formulation}) of the optimization procedure boxplots for 30 runs starting from a random field every run.}
\label{fig:WE_sln_1_1}
\end{figure}

As the exact solution, we take an analytical solution from the Wolfram Mathematica 12.1 software. The solution has the analytical form and is taken at the grid points for every grid used in the optimization process. We record the optimization time and the mean average error (MAE). Time and MAE boxplots have the form Fig.~\ref{fig:WE_sln_1_1}.

The algorithm converges to an analytical solution, and the initial field does not affect the optimization process. Thus, we can say that algorithm converges in this case. Moreover, it shows the stability property since it does not depend on an initial field.

We see that the optimization time is high (more than two hours for a run). The reasons have multiple directions. It could be the optimization method speed, non-optimal initial approximation, and many others. The latter influence of the initial approximation is analyzed below in Sec.~\ref{subsec:initial_field}.

\subsubsection*{Heat equation}

The different types of equations usually require different types of finite-difference schemes and different grid proportions. In contrast with the hyperbolic wave equation, we take the parabolic heat equation boundary problem in the form Eq.~\ref{eq:HE_formulation}.

\begin{equation}
    \begin{array}{cc}
         \frac{\partial u(x,t)}{\partial t}-\frac{\partial^2 u(x,t)}{\partial x^2}=0 \\
         u(-8,t)=u(8,t)= \sin (\frac{\pi}{10} t)\\
         u(x,0)=\sin(\frac{\pi}{8} x) \\
         (x,t) \in [-8,8] \times [0,10]=\Omega
    \end{array}
    \label{eq:HE_formulation}
\end{equation}

\begin{figure}[h!]
\centering
\begin{subfigure}{.5\textwidth}
  \centering
  \includegraphics[width=0.95\linewidth]{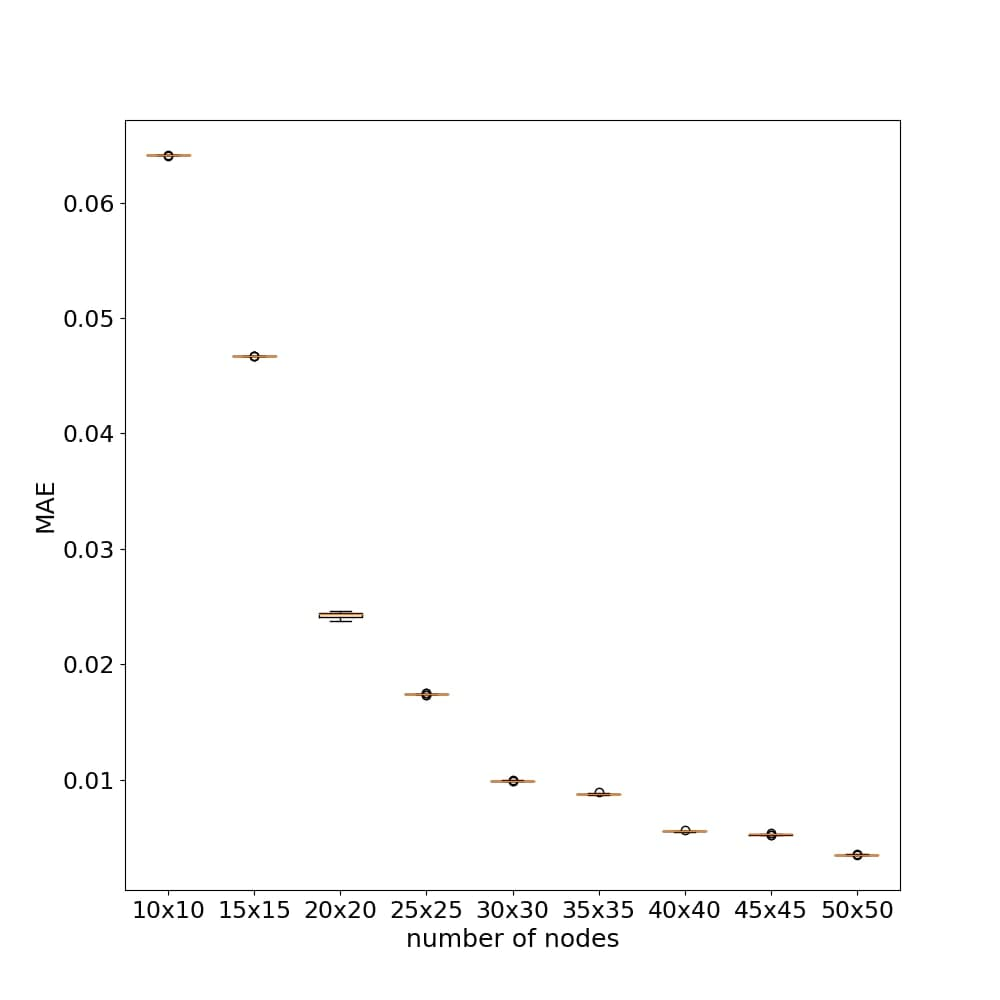}
  \caption{MAE}
\end{subfigure}%
\begin{subfigure}{.5\textwidth}
  \centering
  \includegraphics[width=0.95\linewidth]{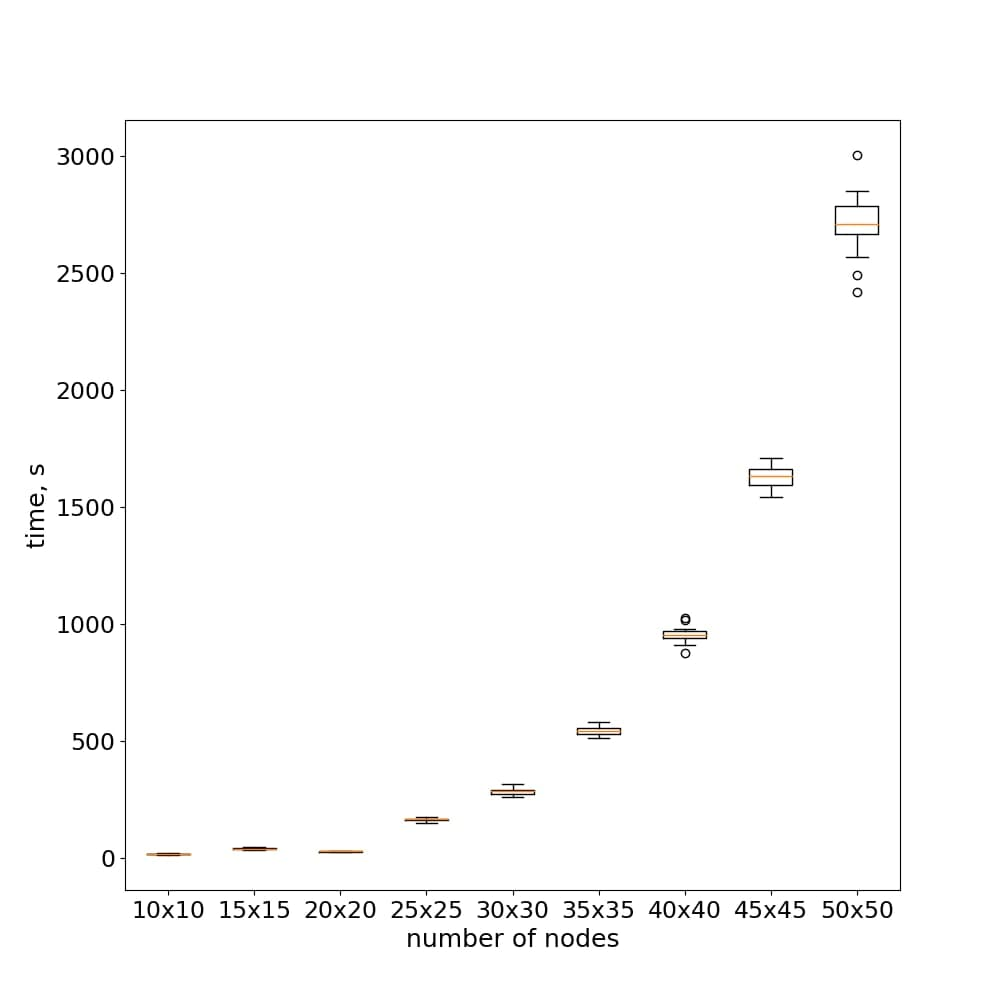}
  \caption{Optimization time}
\end{subfigure}
\caption{Convergency and time complexity for the wave equation (Eq.~\ref{eq:HE_formulation}) of the optimization procedure boxplots for 30 runs starting from a random field every run.}
\label{fig:HE_sln_1_1}
\end{figure}

The same optimization strategy was applied to the heat equation. Every resolution from $10 \times 10$ uniformly distributed points to $50 \times 50$ points with the step of 5 points. We make 30 runs of the algorithm for each number of discretization points starting from the random field. The exact solution was taken from the Wolfram Mathematica 12.1 software. The solution has the analytical form and is taken at the grid points for every grid used in the optimization process. Again we record the optimization time and the mean average error (MAE). Optimization time and MAE boxplots have the form Fig.~\ref{fig:HE_sln_1_1}.

For the heat equation Fig.~\ref{fig:HE_sln_1_1} solution we see the same tendencies as for the Fig.~\ref{fig:WE_sln_1_1}. Additionally, we can conclude that the optimal differentiation realization will speed up the process since the heat equation has less space-derivative order and computes faster.

\subsection{Initial field approximation}
\label{subsec:initial_field}

In this subsection, all experiments are shown only for wave equation Eq.~\ref{eq:WE_formulation} as for the most illustrative case.

For the optimization algorithm, an initial approximation is essential to decrease the optimization time. As the right initial field, one may use the interpolate (or approximate) field obtained from a coarser grid. In experiments, we use several interpolation and approximation algorithms to determine the better to use and then do a final experiment to prove that the approach is working.

As the interpolation algorithms, we take interpn from scipy.interpolate \cite{virtanen2020scipy} package, Rbf class from scipy.interpolate package and a deep neural network.

Interpolation algorithm interpn was taken with default parameters, for Rbf were chosen following parameters: method='linear', smooth=10. 

The neural network (NN) has the following architecture: min-max scaling, two-dimensional input layer with ReLU activation, three dense layers of size 256,64,1024 respectively with ReLU activation, and a one-dimensional output layer with ReLU activation. The architecture was inspired by the multiple runs of the AutoKeras package.

To assess the performance we make 30 runs using following scenario. Starting from random $10 \times 10$ field we obtain $10 \times 10$ solution and interpolate it to $15 \time 15$ grid to obtain initial field for $15 \times 15$ optimization algorithm and continue this process up to $50 \times 50$ resolution. Mean optimization time with 95\% confidence interval is shown in Tab.~\ref{tab:interp_time_1_1}

\begin{table}[h!]
\begin{tabular}{|l|l|l|l|l|l|l|}
\hline
Method & $25 \times25$ & $30 \times30$ & $35 \times35$ & $40 \times40$ & $45 \times45$ & $50 \times 50$ \\ \hline
random& $      0.02 \pm 0.001$& $      0.05 \pm 0.002$& $      0.11 \pm 0.007$& $      0.23 \pm 0.016$& $      0.49 \pm 0.004$& $      0.97 \pm 0.02$ \\ \hline
interpn& $      0.02 \pm 0.001$& $      0.04 \pm 0.001$& \cellcolor[HTML]{009901}$      0.1 \pm 0.001$& $      0.23 \pm 0.004$& \cellcolor[HTML]{009901}$      0.47 \pm 0.007$& \cellcolor[HTML]{009901}$      0.94 \pm 0.01$ \\ \hline
RBF& $      0.02 \pm 0.001$& $      0.04 \pm 0.001$& $      0.11 \pm 0.001$&\cellcolor[HTML]{009901}$      0.22 \pm 0.005$& $      0.48 \pm 0.007$& $      0.96 \pm 0.01$ \\ \hline
NN& $      0.02 \pm 0.002$& $      0.05 \pm 0.003$& $      0.11 \pm 0.013$& $      0.23 \pm 0.010$& $      0.50 \pm 0.035$& $      0.97 \pm 0.03$ \\ \hline
\end{tabular}
\caption{Mean max possible time ratio (1 is the maximum time spent on optimization during the all runs) with 95\% confidence interval (lower the better).}
\label{tab:interp_time_1_1}
\end{table}

We emphasize that random time is shown only in illustrative matters and for maximum possible time computation. From Tab.~\ref{tab:interp_time_1_1} it is seen that scipy interpn provides the best possible interpolation. To simulate the real case scenario, we use interpn to interpolate $25 \times 25$ solution to $50 \times 50$ grid and use it as the initial field for optimization. Resulting mean (30 runs average) optimization time is shown in Tab.~\ref{tab:interp_time_final_1_1}

\begin{table}[h!]
\centering
\begin{tabular}{|l|l|l|l|}
\hline
Interpolation method & $25 \times 25$ & $50\times50$ & Total                        \\ \hline
random               & -     & 7848  & 7848                         \\ \hline
scikit interpn       & 156   & 7008  & \cellcolor[HTML]{009901}7164 \\ \hline
\end{tabular}
\caption{Average (30 runs) time of optimization (sec).}
\label{tab:interp_time_final_1_1}
\end{table}

In this case, we obtain a speed-up of $(7848-7164)/7848 \approx 9 \%$ on average only by changing the initial field. 

Therefore, changing the initial field for optimization may speed-up the process. Nevertheless, we see that the initial field is only a small part of the algorithm's computational complexity.

\subsection{Finite-difference scheme choice}

To show the algorithm convergency for the vast number of cases, we introduce the scheme based on \cite{fornberg1988generation} of fourth-order (approximation order of $O(h^4)$ on a uniformly spaced grid with step $h$) in the form of Eq.~\ref{eq:forward_scheme_second}.

\begin{equation}
   \begin{array}{cc}
    u'_f(x) = -\frac{3}{2} \frac{1}{h} u(x)+2 \frac{1}{h} u(x+h)-\frac{1}{2} \frac{1}{h} u(x+2h)\\
    u'_b(x) = \frac{3}{2} \frac{1}{h} u(x)-2 \frac{1}{h} u(x-h)+\frac{1}{2} \frac{1}{h} u(x-2h)
   \end{array}
   \label{eq:forward_scheme_second}
\end{equation}

For the interior points we use  scheme Eq.~\ref{eq:center_scheme_second}.

\begin{equation}
u'_c(x) =\frac{1}{2}(u'_f(x) + u'_b(x))
   \label{eq:center_scheme_second}
\end{equation}

Increasing of the order of the finite-difference scheme approximation obviously affects the time of the optimization negatively using maximum time ratio from previous section (used in Tab.~\ref{tab:interp_time_1_1}) we obtain time ratios for optimization procedure with scheme Eq.~\ref{eq:forward_scheme_second}-Eq.\ref{eq:center_scheme_second} shown in Tab.~\ref{tab:interp_time_2_2}

\begin{table}[h!]
\begin{tabular}{|l|l|l|l|l|l|l|}
\hline
Method & $25 \times25$ & $30 \times30$ & $35 \times35$ & $40 \times40$ & $45 \times45$ & $50 \times 50$ \\ \hline
random& $      0.03 \pm 0.001$& $      0.06 \pm 0.003$& $      0.13 \pm 0.005$& $      0.28 \pm 0.011$& $      0.60 \pm 0.008$& $      1.14 \pm 0.011$ \\ \hline
interpn& $      0.02 \pm 0.001$& $      0.05 \pm 0.002$& $      0.11 \pm 0.002$& $      0.25 \pm 0.006$& $      0.53 \pm 0.011$& $      1.08 \pm 0.011$ \\ \hline
\end{tabular}
\caption{Time ratio for fourth-order difference scheme (1 is the maximum time spent on optimization during the all runs in Tab.\ref{tab:interp_time_1_1}) with 95\% confidence interval (lower the better).}
\label{tab:interp_time_2_2}
\end{table}

It is seen that using the more computationally complex scheme results in a significant increase of the optimization time. It is known from the classical analysis that most of the classical methods achieve maximum error at the boundaries. Therefore , to reach the optimization time compromise we use second-order scheme Eq.~\ref{eq:center_scheme_first} for interior points and fourth-order scheme Eq.~\ref{eq:forward_scheme_second} for the boundaries. It results with the optimization time shown in Tab.~\ref{tab:interp_time_1_2}.

\begin{table}[h!]
\begin{tabular}{|l|l|l|l|l|l|l|}
\hline
Method & $25 \times25$ & $30 \times30$ & $35 \times35$ & $40 \times40$ & $45 \times45$ & $50 \times 50$ \\ \hline
random& $      0.02 \pm 0.001$& $      0.05 \pm 0.005$& $      0.12 \pm 0.011$& $      0.26 \pm 0.005$& $      0.52 \pm 0.009$& $      1.04 \pm 0.024$ \\ \hline
interpn& $      0.02 \pm 0.001$& $      0.04 \pm 0.002$& $      0.10 \pm 0.003$& $      0.22 \pm 0.003$& $      0.45 \pm 0.005$& $      0.93 \pm 0.005$ \\ \hline
\end{tabular}
\caption{Time ratio for the combined difference scheme (1 is the maximum time spent on optimization during the all runs in Tab.\ref{tab:interp_time_1_1}) with 95\% confidence interval.}
\label{tab:interp_time_1_2}
\end{table}

We note that time is less important than an overall error with respect to the analytical solution. The resulting table of the mean averaged errors (MAE) for all schemes shown above are gathered in Tab.~\ref{tab:interp_MAE}

\begin{table}[h!]
\centering
\begin{tabular}{|l|l|l|l|l|l|l|}
\hline
Method & SO & BO & $30 \times30$  & $40 \times40$ &  $50 \times 50$ \\ \hline
random& 2 & 2& $0.0069 $& $0.0054 $& $0.0050 $ \\ \hline
interpn& 2 & 2& $0.0069 $& $0.0038 $& $0.0024 $ \\ \hline
RBF& 2 & 2& $0.0069 $& $0.0038 $& $0.0024 $ \\ \hline
NN& 2 & 2& $0.0069 $& $0.0046 $& $0.0051 $ \\ \hline
random& 2 & 4& \cellcolor[HTML]{009901}$0.0036$& $0.0044 $& $0.0040 $ \\ \hline
interpn& 2 & 4& \cellcolor[HTML]{009901}$0.0036$& \cellcolor[HTML]{009901}$0.002$& \cellcolor[HTML]{009901}$0.0012 $ \\ \hline
random& 4 & 4& $0.007 $& $0.006$& $0.0062 $ \\ \hline
interpn& 4 & 4& $0.007 $& $0.0048 $& $0.0041 $ \\ \hline
\end{tabular}
\caption{Mean averaged error for all runs (SO - order of the scheme for interior points, BO - order of the scheme at the boundary points)(lower the better).}
\label{tab:interp_MAE}
\end{table}

It is seen from Tab.~\ref{tab:interp_time_1_2} and Tab.~\ref{tab:interp_MAE} that the boundary points scheme approximation gives an insignificant time increase in case of random initial field and even time boost combined with the initial field interpolation. More important that it gives the lowest MAE with respect to the analytical solution. It entirely agrees with the classical analysis conclusion - boundary points approximation is crucial for the numerical methods.

\section{Conclusions}
\label{sec:conclusions}

In the paper, we propose a numerical method of PDE solution based on optimization methods. It has the following advantages:

\begin{itemize}
    \item It can solve PDEs without the involvement of an expert, which is most useful for data-driven methods.
    \item It has good precision every for in experimental realization
    \item It can be easily parallelized
    \item It has a flexible modular structure. The modules could be replaced to achieve the better speed or better precision
\end{itemize}

 It is also seen that the optimization time is the main drawback of the method's experimental realization. We propose the following optimization speed-up directions:

\begin{itemize}
    \item usage the power of GPU to make the optimization using fast memory and built-in matrix instructions
    \item Better usage of initial approximation
    \item More intelligent use of the numerical differentiation
    \item Usage of matrix-based optimization methods
\end{itemize}

All experimental data and script that allows reproducing experiments are available at the GitHub repository \footnote{github.com/ITMO-NSS-team/FEDOT.Algs/tree/master/PDE\_solver}.

%
%
%
 \bibliographystyle{splncs04}
\bibliography{mybibliography}

\end{document}